\documentclass[12pt]{article}

\usepackage[intlimits]{amsmath} 
\usepackage{amsthm,mathrsfs}    
\usepackage[comma,authoryear]{natbib}
\usepackage{url}
\numberwithin{equation}{section}
\usepackage{amssymb,amsmath}
\usepackage{subfigure}
\usepackage{bm}
\usepackage{booktabs} 
\usepackage{graphicx,graphics,epsfig}
\usepackage{float}

\theoremstyle{plain}
\newtheorem{thm}{Theorem}

\newcommand{\bthm}{\begin{thm}}
\newcommand{\ethm}{\end{thm}}

\newcommand{\bpf}{\begin{proof}}
\newcommand{\epf}{\end{proof}}

\theoremstyle{definition}

\newcommand{\opn}{\operatorname}
\newcommand{\KLP}{\operatorname{KLP}}

\usepackage[a4paper,text={16.5cm,25cm},centering]{geometry}
\usepackage[compact]{titlesec}
\usepackage{deep-macro}

\begin{document}
\begin{center}
{\Large {\bf LP Mixed Data Science : Outline of Theory }}
\\[.2in]
 Emanuel Parzen$^1$, Subhadeep Mukhopadhyay$^2$,\\
$^1$Texas A\&M University, College Station, TX, USA\\
$^2$Temple University, Philadelphia, PA, USA\\
~October, 2013 \\[.3in]
\end{center}
\vskip2em

\medskip\hrule height 1pt

\tableofcontents

\vspace*{.6in}

\newpage

\setcounter{section}{-1}
\section{Introduction}
Statistics journals have great difficulty accepting papers unlike those previously published.  For statisticians with  new big ideas a practical strategy   is  to publish them in  many small applied studies which enables one to provide references to work of others.  This essay outlines the  many concepts,  new theory, and  important algorithms of our new culture of statistical science called LP MIXED DATA SCIENCE.   It provides  comprehensive solutions to  problems of  data analysis and nonparametric modeling of   many   variables that are continuous or discrete, which  does not yet have a large literature.  It develops a new modeling approach to nonparametric estimation of the multivariate copula density.  We discuss the theory which we believe is very elegant (and can provide a framework for  United Statistical Algorithms, for traditional Small Data methods and Big Data methods). The utility of the theory will be demonstrated elsewhere by series of real applications.
\section{LP Methods, Mixed Variables (X,Y), Mid-Distribution}
By mixed variables we mean $X,Y$ can be continuous or discrete. We will denote distribution functions by $F(x;X),  F(y;Y)$ and Quantile functions by $Q(u;X), Q(v;Y), 0<u<1, 0<v<1$.

One important mixed data  problem is CLASSIFICATION: $Y$ binary $0,1$; $X$ continuous. The goal: Estimate nonparametrically $\Pr[Y=1|X=x=Q(u;X)] =\Pr[Y=1|F(X;X)=u]$ as a function of  values $u$ of rank transform $F(X;X)$.

LP METHODS: {\bf L} stands for extension of L rank statistics, L moments.{\bf P} stands for Parzen, quantiles, mid-distribution, comparison densities ,  orthonormal score functions $T_j(x;X)$ custom built for each $X$ as functions of mid-ranks. Mid-distribution is defined as $\Fm(x;X)=F(x;X)-.5 p(x;X)$, probability mass function $p(x;X)=\Pr[X=x]$, notation introduced by Parzen (1960)``Modern Probability Theory and its Applications''. Sample mid-distribution from sample of size $n$ computed in \texttt{R} by
$\tilde{F}^{\mbox{mid}}(X;X)=(\mbox{Rank}(X)-.5)/n$.

LOOKING AT  UNIVARIATE DATA: Quantile plot of sample with distinct values $x_j$  is scatter diagram $(u_j=\Fm(x_j;X), x_j=Q(u_j;X))$. We define Mid-quantile $\Qm(u;X), 0<u<1$, linearly connect $(u_j,x_j)$. Normal Q-Q plot is scatter diagram  $(Q(u_j;Z),x_j)$,  where $Z \,\mbox{Normal}(0,1)$.

LOOKING AT BIVARIATE DATA $(X,Y)$: We recommend to display three scatter plots $(X,Y)$, $(\Fm(X;X), Y)$, $(\Fm(X;X),\Fm(Y;Y))$.

COMPARISON DENSITY: To test if two continuous distributions $F(x)$ and $G(x)$ are equal we compare \emph{not} their difference but their ratio by defining comparison distribution $D(u;G,F)=F(Q(u;G))$ with comparison density $d(u;G,F)=f(Q(u;G))/g(Q(u;G))$, likelihood ratio evaluated at $x=Q(u;G), u=G(x)$.

\section{Score Function Construction, First Step of Modeling and Algorithm}
For a distribution $F$ we construct orthonormal score functions $T_j(x;F)$ as functions of $x$; denoted $T_j(X;X)$ as transforms of $X$.
Define  $T_1(X;X)=T_1(X;F)=(\Fm(X;X)-.5)/\si(\Fm(X;X)),$
\begin{thm}[Parzen (2004)]
$\Ex[\Fm(X;X)]=.5,$ and  $\Var[\Fm(X;X)]=(1/12)(1-\sum_x p^3(x;X))$
\end{thm}

For Y 0,1 , $p=\Pr[Y=1]$, and $\mbox{Var}[\Fm(X;X)]=pq/4=(1/12)(1-p^3-q^3)$.

Define Score functions $T_j(X;X)=T_j(X;F)$ constructed by  Gram Schmidt orthonormalization of powers of $T_1(X;X)$. Define  Score functions $S_j(u;X)=T_j(Q(u;X);X)$.

\section{Legendre Polynomial Score Functions for $X$ Continuous}
Orthonormal  Legendre polynomials $\Leg_j(u)$ on interval $0<u<1$ are constructed by Gram Schmidt orthonormalization of $1, u, u^2, \cdots .$ First four Legendre polynonomials are given by:
\beas
\Leg_0(u)&=&1 \\
\Leg_1(u)&=& \sqrt{12} (u-.5)\\
\Leg_1(u)&=& \sqrt{5} (6u^2-6u+1)\\
\Leg_3(u)&=&\sqrt{7} (20u^3-30u^2+12u-1)\\
\Leg_4(u)&=&3(70u^4-140u^3+90u^2-20u+1)
\eeas
For $X$ continuous our score functions are $S_j(u;X)=\Leg_j(u)$, $T_j(x;X)=\Leg_j(F(x;X))$.

We argue that we can model $\Pr[Y=1|X=x]$ or $\mbox{logodds} \Pr[Y=1|X=x]$ not as a linear function of powers of $x$ but  as a linear function of  orthonormal score functions of $\Fm(x;X)$ which yields model of $\mbox{logodds} \Pr[Y=1|X=Q(u;X)]$ as linear function of orthonormal score functions $S_j(u;X)$.  We identify score functions \emph{before} parameter estimation (which comes first, parameters or sufficient statistics?). Logistic regression algorithms can be applied to numerically compute parameters of identified score functions.

\section{Mixed Statistical Data Science, Unify Small and Big Data}
Data scientists dispute  claim by statisticians that data science is just a sexier name for applied statistical science.  They agree with  Nate Silver (2013 JSM)  that \emph{the ultimate goal is quality applied  research that gets practical answers}. Our view is that data science has many cultures, including applied statistical science, machine learning, statistical learning.  \emph{It is useful to distinguish} between those who  strive for utility, and those who aim for utility and elegance  (applying RKHS (reproducing kernel Hilbert spaces), regularization, and sweep regression).  Instruction in data science  often  presents \emph{not} unified theory  but  algorithms (formulas for final answers) to be imitated rather than understood, including  ``looking at the dat''.

CHALLENGE! To understand differences between cultures of data science study their approach to two sample inference, classification, for data ($Y$ binary $0-1$, $X$ continuous).

A VISION FOR FUTURE OF STATISTICS:  We propose that applied statistical science  can be  more scientific (and less art) when it emphasizes:

{\bf (1)} unified methods and graphical analysis  that work for small data and big data (analogies between analogies), and

{\bf (2)}  awareness of the history and scope of statistical methods (confirmatory and exploratory), as  outlined below.

\vskip.5em

{\bf A.	Modern Probability:}
Think frequentist, compute axiomatically (Kolmogorov 1933)
Awareness of definition  mixed conditional probability $\Pr[Y=y|X=x]$ for $Y$ discrete, $X$ continuous.
Parzen (1960) teaches modern probability \emph{without} measure theory. NEW! Mid-probability theory: mid-distribution inversion, convergence.
\vskip.5em

{\bf B.	Parametric Inference (Objective):}
Think Bayesian (parameter probability), compute frequentist confidence quantiles to combine Fisher and Neyman.

\vskip.5em
{\bf C. Nonparametric Inference (Quantiles, Ranks):}
Think nonparametrically, compute parametrically by models selected nonparametrically using information criteria.

\vskip.5em
{\bf D.	Statistical Data Science, Unify Cultures of Small and Big Data:}
COMPRESSION!   goal of  high dimensional data analysis, ``statistics is like art, like dynamite, the goal is compression'';  reduce  number of  influential variables,  number of sufficient statistics that summarize  massive
Data.

\section{Quantile Mechanics}
The quantile function $Q(u;X)$ of a random variable $X$ is defined $Q(u;X)=\inf \{x: F(x;X) \geq u\},$ $0<u<1.$
Call $u$ probable if there exists $x$ such that $F(x;X)=u$. Verify $Q(u;X)$ equals  $x(u)$,  smallest $x$ such  that $F(x;X)=u$. Verify that the exact inverse property $Q(F(x;X);X)=x$ holds for $x=x(u)$. With probability one, an observed value of $X$ equals  $x(u)$ for some probable $u$.

\begin{thm}
Any random variable $X$ has the property that with  probability $1$,  Q(F(X;X);X)=X. With probability $1$, a function $h(X)$ of $X$ equals a function of rank transform $F(X;X)$ since $h(X)=hQ(F(X;X))$, defining $hQ(u)=h(Q(u;X))$.
\end{thm}
We now state a REMARKABLE fundamental theorem on conditional expectation by rank transform.

\begin{thm}
With probability $1$,  $\Ex[Y|X]\,=\,\Ex[Y|F(X;X)]\,=\,\Ex[Y|\Fm(X;X)].$
\end{thm}
This theorem has direct relevance to NONPARAMETRIC REGRESSION given by the following theorem.
\begin{thm}
$\Ex[Y|X=x=Q(u;X)]$ can be approximated by a linear combination of orthonormal score functions $T_j(X;X)$ with coefficients $\Ex[Y T_j(X;X)] =\LP(j,0;X,Y)$.
\end{thm}

\begin{thm}[LP Representation]
$E[g(Y)|X]=E[g(Y)|F(X;X)]=\sum_j T_j(X;X)  \Ex[g(Y) T_j(X;X)]$.
\end{thm}

\begin{thm}
$\Var(\Ex[Y|X])\,=\,\sum_{j>0} |\LP(j,0;X,Y)|^2$.
\end{thm}

We present below  LP representation of $\Var(X)$ equivalent to  letting $Y=X$. We explore below consequences of LP representation of variance and concept of tail index of the distribution of a random variable $X$.

\begin{thm}[Parzen (1979)]
$Q(u;g(X))=g(Q(u;X))$ for $g$ quantile-like  (non-decreasing and left continuous function).
\end{thm}

\begin{thm}
Let U denote Uniform(0,1) variable. In distribution X=Q(U;X).  For X continuous,  In distribution F(X;X)=U because  F(Q(u;X);X)=u, all u. $X$ discrete $F(Q(u;X);X)=u$ for probable $u$.
\end{thm}

\begin{thm}[Conditional Quantile]
Because $Y=Q(F(Y;Y);Y)$ with probability $1$  the conditional quantile function $Q(u;Y|X)$ can be computed $Q(u;Y|X)=Q(Q(u;F(Y;Y)|X);Y)$ by first computing the  conditional given X quantile of the rank transform F(Y;Y).
\end{thm}
We simulate this by estimating the comparison density $d(v;Y,Y|X)$, a very important concept defined below.

\section{LP Moments of $X$, Tail Index}
Quantile formulas  every statistician should know!   Mean $\Ex[X]=\Ex[Q(U;X)]$, Variance $\Var[X]=\Var[Q(U;X)]$.
Other measures of location and scale can be defined in terms of quartiles Q1, Q3, and median Q2. Mid quartile  MQ=.5(Q1+Q3), quartile deviation DQ=2(Q3-Q1)  approximate slope at Q2. Next we introduce the concept of informative quantile function, a powerful exploratory data analysis tool. Define INFORMATIVE QUANTILE $\QI(u;X)$, quantile of $\QI(X)=(X-\mbox{MQ})/\mbox{DQ}$. An observed value $X$ is called Tukey outlier if $|\QI(X)|>1$.

GINI  Method: Measure of scale of $X$ continuous is $\Ex[X(F(X;X)-.5)]= (1/4)\Ex[|X-X'|]$, where $X$ and $X'$ independent identically distributed. It has  a long history of theory and application  under name Gini coefficient. Definition for $X$ discrete is \emph{not} obvious; most accepted answer is equivalent to our general definition $\LP(1;X)=\Ex[X(\Fm(X;X)-.5)]/\si(\Fm(X;X))$.

LP MOMENTS:  $\LP(j;X)$ of $X$:  $X$ continuous, $\LP(j;X)=\Ex[X  \Leg_j(F(X;X))]$. All $X$ define  $\LP(j;X)=\E[X T_j(X;X)]$.
For a distribution $F$ with quantile $Q$, define LP moments $\LP(j;F)=\Ex[Q(U;F) T_j(Q(U;F);F)]$. Note for $X$ continuous,  $T_j(Q(U;F);F)=\Leg_j(U)$.

\begin{thm}[LP Representation of Quantile Function]
$Q(u;X)\,=\, \sum_j S_j(u;X) \LP(j;X)$.
\end{thm}

LP Representation of Variance: $\Var[X]=\sum_{j>0} |\LP(j;X)|^2$.

An empirical representation or estimator from data of $Q(u;X)$ when $X$ is continuous  is
$\widehat Q(u;X)= \sum_{\mbox{selected}\, j}  \Leg_j(u) \widehat\LP(j;X).$

\begin{thm}[Orthonormal Representation of random variable X]
With probability $1$, $X=\sum_j T_j(X;X) \LP(j;X)$, and $\cZ(X)=\sum_j T_j(X;X) \LP(j;\cZ(X))$.
\end{thm}

TAIL BEHAVIOR OF DISTRIBUTIONS: Parzen (1979) classifies tails of distributions into  short, long, medium (medium-short, medium-medium, medium-long). Normal is medium-short.  A statistical joke: the tails (ends) justify the means (location estimator).

LP SHORT TAIL DISTRIBUTIONS. Recall $\cZ(X)=(X-\Ex[X])/\si(X)$. LP  Tail index of $X$ is smallest $m$ that
$\sum_{0< j \leq m} |\LP(j;\cZ(X))|^2 > .95.$

Threshold $.95$ is chosen because Normal barely satisfies it. One could use  threshold $.99$ to choose number of LP moments in an empirical representation of $Q(u;X)$.

\begin{thm} For $Z$ Normal(0,1),
$\LP(1;Z)=\sqrt{12} \Ex[Q(U;Z) U]=\sqrt{12}\Ex[fQ(u;Z)] =\sqrt{3/\pi}=.977$,
Density quantile function $fQ(u;Z)=f(Q(u;Z);Z)$.  Verify  score  $J(u;Z)=-(fQ(u;Z))'=Q(u;Z)$.
\end{thm}

\begin{thm}
Normal $X$ is short tailed. Uniform $X$ has $\LP(1;\cZ(X))=1$.
\end{thm}

Goodness of fit test of Normality based  on  sample $\LP(1;\cZ(X))$ is analogous to Shapiro Wilk  test $\Ex[Q(\Fm(X;X);Z) \cZ(X)]= \Ex[Z(X) \mbox{Hermite}_1(\Fm(X;X))]$.

L MOMENTS: Our concept of LP moments extends to  discrete variables concept of L moment advocated by Hosking (1990)

LP CRITERION: Identify monotonic transform $g(X)$ which is short tailed, i.e., satisfies $E[\cZ(g(X)) \cZ(\Fm(X;X))] > .975.$

\section{LP Comoments of (X,Y), Covariance Matrix of Score functions}

Define  $\LP(j,k;X,Y)=\Ex[T_j(X;X) T_k(Y,Y)]$.

Our concept LP comoment extends concept L comoment introduced by Serfling and Xiao (2007). For analogies with multivariate analysis use Covariance matrix $\KLP(X,Y)$ of vectors ${T_j(X;X)}$ and ${T_k(Y;Y)}$, selected $j,k$. Note $\KLP(X;X)$ and $\KLP(Y;Y)$ are identity.
Our criterion $\LPINFOR(X,Y)$ for independence of $X$ and $Y$, correspondence analysis, canonnical correlations squared are  based on  eigenvalues of LP-Coherence(X,Y)=$\KLP(X,Y) \KLP(Y,X)$.

Define $\LPINFOR(X,Y)= \sum_{j,k} |\LP(j,k;X,Y)|^2 =\mbox{trace}\, \mbox{LP-Coherence}(X,Y)$.

\begin{thm}[Correlation Representation of LP Comoments]
From orthonormal representations for $\cZ(X), \cZ(Y)$ obtain
$R(X,Y)=\sum_{j,k>0} \LP(j;\cZ(X))\,\LP(j,k;X,Y)\, \LP(k;\cZ(Y))$, and
$R^2(X,Y) \leq \sum _{j,k} |\LP(j,k;X,Y)|^2=\LPINFOR(X,Y)$.
\end{thm}

$\LPINFOR(X,Y)$ defined above is an information measure of dependence.

\begin{thm}[REMARKABLE approximate equality of  Pearson and Spearman Correlation]
For $X,Y$ short tailed approximately $R(X,Y)=\LP(1;X)\,\LP(1,1;X,Y)\,\LP(1;Y)$.
This can be verified directly for $(X,Y)$ bivariate normal.
Spearman Correlation can be defined for X,Y mixed $\operatorname{RSPEARMAN}(X,Y)=\LP(1,1;X,Y)$.
\end{thm}
Spearman correlation equals Pearson correlation of mid-rank transform   $\Fm(X;X)$ and $\Fm(Y;Y)$.

TIES IN DATA: Estimation of Spearman correlation is difficult  when data has many ties. Our definition of sample Spearman correlation avoids this problem.

\begin{thm}
$R(X,Y)=\operatorname{RSPEARMAN}(X,Y)$  for X, Y both uniform or both binary 0,1.
\end{thm}
LPINFOR, CHISQUARED  EMPIRICAL INFORMATION STATISTIC FOR INDEPENDENCE: We define data-driven nonlinear measure of dependence by
$\LPINFOR(X,Y)=$ $\sum_{\mbox{significant}\, j,k} |\LP(j,k;X,Y)|^2$.

CHISQUARE TEST FOR INDEPENDENCE: Also call LPINFOR  Chi-square divergence to test independence of $X$ and $Y$, extension of  usual chi-squared statistic which we express below as integral of square of discrete bivariate copula density function.

\section{Skew-G Distribution, Goodness-of-fit}
To estimate probability density $f(x;X)=F'(x;X)$ of continuous $X$  a popular method is kernel density estimation. A powerful alternative method  starts with a  parametric model $G(x)$ with density $g(x)$.

A model  for unknown $f(x;X)$, called Skew-G model, is $f(x;X)=g(x) d(G(x))$, where $d(u)$ is the comparison density $d(u)=d(u;G,F(.;X)) =f(Q(u;G);X)/g(Q(u;G))$. Probability density of  G-transform $G(X)$ with distribution function $F(u;G(X))=D(u)=D(u;G,F(.;X))=F(Q(u;G);X)$,
called comparison distribution.

Estimator $\widehat d$  of  $d$, which provides estimator $\widehat f$ of $f$, can be formed by orthonormal score function representation or by exponential (maximum entropy)  model. Verify $d(u)-1= \sum_j \Leg_j(u)\, \langle d, \Leg_j \rangle$,
where $\langle d, \Leg_j \rangle$ is inner product of density $d$ and $\Leg_j(u)$ which is evaluated as the sample mean $\Ex[\Leg_j(G(X))]$. An empirical estimator of $d(u)$ has representation $\widehat d(u) -1= \sum_{\mbox{selected}\, j} \Leg_j(u) \Ex[\Leg_j(G(X))]$.

COMPONENT TESTS GOODNESS OF FIT of continuous $G(x)$ to $F(x;X)$ is tested by values of $\Ex[\Leg_j(G(X))]$ which we call component tests.

Define G COMPONENTS: $\opn{Comp}_j(X;G)=\Ex[T_j(X;G)]$, where expectation uses (sample) distribution of $X$. Note $\opn{Comp}_j(X;G)=0$ when $X$ has true distribution $F$ equal to $G$.

\begin{thm}[Component Representation of Comparison Density] We have the following orthonormal expansion:
$d(u;G,F)\,=\,\sum_j  T_j[Q(u;G);G]  \opn{Comp}_j(X;G).$
\end{thm}

DISCRETE MODELING AND GOODNESS OF FIT: When  we observe a sample of discrete $X$ with probable values $x$, a null hypothesis is probability mass function $g(x)$, usually denoted $p_0(x)$. Sample probability mass function is $\widetilde p(x)$. Sample comparison density is $\widetilde d(u)=d(u;G,\widetilde F)=\widetilde p(Q(u;G))/g(Q(u);G)$.

First component statistic is $\widetilde{\Ex}[\cZ(\Gm(X))]$. First estimator of true $p(x)$ is $\widehat{p}(x)=g(x)\Big\{1+\cZ(\Gm(x)) \widetilde{\Ex}[\cZ(\Gm(X))] \Big\}.$

EXAMPLE: ED JAYNE'S DIE: Estimate probabilities of $6$ sided die when observed sample mean of $X$ equals $4.5$; a fair die has population mean $3.5$. Note we are \emph{not} given value of $n$, sample size.

EXAMPLE SPARSE CHI-SQUARED  (large $p$, small $n$). Let discrete $X$ have $p=20$ outcomes, Sample probabilities are $.75, .25$ for first two outcomes. Model for population model probabilities are $.25$ for first two outcomes and $1/36$ other $18$ outcomes. Sample size is $n=20$. Test model and estimate from data observed true probabilities of outcomes using model as parametric start.

\section{Copula Density, Conditional Comparison Density, Comparison Probability Bayes Rule}
When $X$ and $Y$ are both continuous, or both discrete,  their joint probability is described  by joint probability density  $f(x,y;X,Y)$ or by joint probability mass function $p(x,y;X,Y)$. When $Y$ is discrete and $X$ is continuous joint probability is described by either side of identity $\Pr[Y=y|X=x] f(x;X)=f(x;X|Y=y) \Pr[Y=y]$ which we call PRE-BAYES THEOREM.

BAYES  RULE: $\Pr[Y=y|X=x]/\Pr[Y=y]=f(x;X|Y=y)/f(x;X).$

COMPARISON PROBABILITY: Define left side of Bayes rule  to be $\Comp[Y=y|X=x]$. Define right side of Bayes rule  to be $\Comp[X=x|Y=y]$.

\begin{thm}[Bayes Rule for MIXED  X,Y]
Bayes Rule using Comparison Probability $\Comp[Y=y|X=x]=\Comp[X=x|Y=y]$.
\end{thm}
	
CONDITIONAL COMPARISON  DENSITY: Let $x=Q(u;X), y=Q(v;Y)$. Define conditional comparison   density  of $X$ given $Y$
$d(u;X,X|Y=Q(v;Y))=\Comp[X=Q(u;X)|Y=Q(v;Y)]$.

Conditional comparison density of $Y$ given $X$ $d[v;Y,Y|X=Q(u;X)]=\Comp[Y=Q(v;Y)|X=Q(u;X)]$.

COPULA DENSITY: $\cop(u,v;X,Y)$ to be common value of above conditional comparison densities.  When $X,Y$ jointly continuous copula density function equals $\cop(u,v;X,Y)= f(Q(u;X),Q(v;Y);X,Y)/f(Q(u;X);X)f(Q(v;Y);Y)$.

Copula distribution $\Cop(u,v;X,Y)=F(Q(u;X),Q(v;Y);X,Y)$.

MODELING (X,Y): Estimate univariate  marginal of $X$, univariate marginal of $Y$, joint copula density of $(X,T)$

\begin{thm}[LP Representation of Copula Density]
$$\cop(u,v;X,Y)-1= \sum_{j,k>0}  \LP(j,k;X.Y)\, S_j(u;X)\, S_k(v;Y),~ 0<u,v<1.$$
Equivalently,
$$\int_{[0,1]^2} \dd u \dd v d(v;Y,Y|X=Q(u;X)) S_j(u;X) S_k(v;Y)\,=\,\LP(j,k;X,Y).$$
\end{thm}

A proof of copula LP representation is provided by representations of conditional copula density and conditional expectations:
\beas d(v;Y,Y|X=Q(u;X) &=& \sum_k   S_k(v;Y)\, \Ex[T_k(Y;Y)|X=Q(u;X)]\\
\Ex[T_k(Y;Y)|X=Q(u;X)]&=& \sum_j   S_j(u;X) \Ex[T_j(X;X) T_k(Y;Y)] \eeas

MODELING DEPENDENCE IN PRACTICE:  The LP representation of the copula density provides data driven estimators of the copula density after constructing custom built score functions and LP comoments.

SLICE PLOTTING OF COPULA DENSITY: Plot for selected values of  $u$, as function of $v$, $\cop(u,v;X,Y)=d(v;Y,Y|X=Q(u;X))$.

CONDITIONAL QUANTILE $Q(v;Y,Y|X=Q(u;X))$ can be simulated from conditional comparison density $d(v;Y,Y|X=Q(u;X))$.

\section{Two Sample Inference, Unify Small and Big Data Modeling, Classification}
Two sample inference is equivalent to ($X$ continuous, $Y$ binary 0-1).  A complete analysis estimates conditional comparison density $d(u;X\, \mbox{pooled sample},\, X|Y=1\,)$. Dependence of  $X$ and $Y$ is measured by $\LPINFOR(X,Y)$.

A quick measure of  independence, equivalent to Wilcoxon linear rank statistic  and Spearman correlation,  is $\LP(1,1;X,Y)=R(\Fm(X;X), \ind\{Y=1\})$, equal to $\Ex[\cZ(\Fm(X;X \,\mbox{pooled})|Y=1, X\,\mbox{in sample}\, 1] \sqrt{\mbox{odds} \Pr[Y=1]}.$ This statistic is asymptotically  $\cN(0,1/n)$ under null hypothesis of independence of $X$ and $Y$.

Traditional two sample Student $t$ statistic to test equality of means $\Ex[X|Y=1]=\Ex[X|Y=0]$ assuming equality of variances $\Var[X|Y=1]=\Var[X|Y=0]$ is equivalent to $T=R(X,Y)/\sqrt{1-R^2(X,Y)}$. $R(X,Y)=\Ex[\cZ(X)\cZ(Y)]=\Ex[\cZ(X)|Y=1] \sqrt{\mbox{odds} \Pr[Y=1]}$.

LPINFOR COMPRESSION: When there are many features $X_m$  one wants to identify a small number of features to use to predict (classify)  the value of  $Y$.  For each feature $X_m$ estimate $\LPINFOR(X_m,Y)$. By plotting ranked  values $\LPINFOR(X_m,Y)$ one can start the process of  identifying the features $X_m$ which are most predictive of $Y$.

\vskip3.3em

\end{document}